\documentclass{amsart}

\usepackage{graphicx,amsfonts, amsmath, amssymb, graphics, setspace, tabu}
\DeclareMathOperator{\R}{Re}
\DeclareMathOperator{\I}{Im}
\DeclareMathOperator{\Arg}{Arg}

\newtheorem{theorem}{Theorem}[section]
\newtheorem{lemma}[theorem]{Lemma}
\newtheorem{corollary}[theorem]{Corollary}
\newtheorem{proposition}[theorem]{Proposition}

\theoremstyle{definition}
\newtheorem{definition}[theorem]{Definition}

\theoremstyle{remark}
\newtheorem{remark}[theorem]{Remark}

\numberwithin{equation}{section}

\begin{document}

\title{$\tau$-Norm-Perfect and $\tau$-Perfect Eisenstein Integers for $\tau=\omega+2$ and $2$}

\author{Carlos Rojas Mena}

\curraddr{Department of Mathematics,
University of Rochester, Rochester, New York 14627}
\email{crojasmena1@gmail.com}

\thanks{This work was funded by the SUAMI program at Carnegie Mellon University, which is supported by the NSA}



\begin{abstract}
Using Robert Spira's \cite{D} definitions of complex Mersenne numbers and the complex sum-of-divisors function, we characterize $(\omega+2)$-norm-perfect and $(\omega+2)$-perfect numbers that are divisble by $\omega+2$ and prove the nonexistence of $2$-norm-perfect numbers that are divisible by $2$ in the Eisenstein integers.
\end{abstract}

\maketitle

\section{Introduction}

Let $\sigma: \mathbb{Z}\rightarrow \mathbb{N}$ be the function defined by the equation

\begin{equation}
    \sigma(n) = \sum_{d|n} d
\end{equation} This function is called the sum-of-divisors function. 

In the integers, a $k$-perfect number is a positive integer $n$ satisfying the equation

\begin{equation}
    \sigma(n)=kn
\end{equation}

The most widely studied $k$-perfect numbers are the 2-perfect numbers which are most commonly known by the name of perfect numbers. The first seven 2-perfect numbers are: $6=1+2+3$, $28=1+2+4+7+14$, $496$, $8128$, $2^{12}(2^{13}-1)$, $2^{16}(2^{17}-1)$, and $2^{18}(2^{19}-1)$. As of today, the mathematical community knows exactly 49 2-perfect numbers in the integers. The largest one has 44677235 digits.

The study of perfect numbers dates as far back as Euclid, who circa 300 B.C, proved that, for primes $p$ such that $2^p-1$ is also prime, the numbers of the form 

\begin{equation}
2^{p-1}(2^p-1)
\label{eq:evenform}
\end{equation} are 2-perfect. Numbers of the form $2^p-1$ are now known as Mersenne numbers. In particular, if $2^p-1$ is  prime, it is called a Mersenne prime.

Around two millennia after Euclid's proof, Euler proved that all even 2-perfect numbers were of the form (\ref{eq:evenform}), thereby characterizing all even 2-perfect numbers in the integers.

\begin{theorem}[Euclid-Euler Theorem] The positive integer $n$ is an even 2-perfect number if and only if $n=2^{p-1}(2^p-1)$ where $2^p-1$ is prime.
\end{theorem}

The purpose of this paper is to characterize all $(\omega+2)$-perfect numbers divisible by $\omega+2$ and all  $(\omega+2)$-norm-perfect numbers divisible by $\omega+2$ in the Eisenstein integers, and to show that there exist no $2$-norm-perfect Eisenstein integers divisible by $2$. We follow Wayne McDaniel's \cite{A} and Kieran Smallbone's \cite{B} approach who provided partial characterizations of $(i+1)$-norm-perfect and $(i+1)$-perfect numbers in the Gaussian integers, and $2$-perfect numbers in the Eisenstein integers, respectively.

This paper is structured as follows. In section \ref{background}, we provide some technical background. In section \ref{perfect}, we present our results on $(\omega+2)$-norm-perfect and $(\omega+2)$-perfect numbers. In section \ref{smallbonesresults}, we prove the nonexistence of $2$-norm-perfect that are divisible by $2$ in the Eisenstein integers that are divisible by $2$. In section \ref{discussion}, we discuss some of the unanswered questions about $\tau-$perfect numbers in quadratic integer rings such that the Gaussian and the Eisenstein.

\section{Background}\label{background}

\begin{definition}[Eisenstein integers] The set $\mathbb{Z}[\omega]=\{a+b \omega\mid a,b \in \mathbb{Z}[\omega]\}$, under the usual operations of addition and multiplication of complex numbers, is the ring of Eisenstein integers, where $\omega=e^{\frac{2\pi i}{3}}=\frac{-1 + \sqrt{-3}}{2}$.
\end{definition}

\begin{figure}
\includegraphics{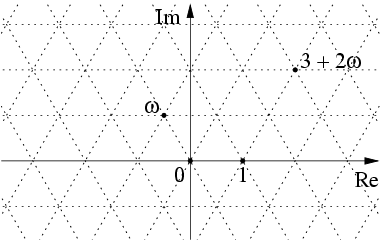}
\caption{The Eisenstein integers form a triangular lattice over the complex plane}
\label{fig:lattice}
\end{figure}

Throughout this paper, it might be helpful for the reader to visualize the Eisenstein integers as a subset of the complex plane. See figure \ref{fig:lattice}. Like the complex plane is partitioned symmetrically into four quadrants, the Eisenstein integers is symmetrically and radially partitioned into six sextants. Each sextant is defined as follows. 

\begin{enumerate}
\item First sextant: $\{\eta \in \mathbb{Z}[\omega] \mid 0\leq \Arg(\eta)<\frac{\pi}{3}\}$
\item Second sextant: $\{\eta \in \mathbb{Z}[\omega] \mid \frac{\pi}{3} \leq \Arg(\eta)<\frac{2\pi}{3}\}$
\item Third sextant: $\{\eta \in \mathbb{Z}[\omega] \mid \frac{2\pi}{3} \leq \Arg(\eta)<\pi\}$
\item Fourth sextant: $\{\eta \in \mathbb{Z}[\omega] \mid  -\pi< \Arg(\eta)<-\frac{2\pi}{3} \text{ or } \Arg(\eta)=\pi\}$
\item Fifth sextant: $\{\eta \in \mathbb{Z}[\omega] \mid  -\frac{2\pi}{3}\leq \Arg(\eta)<-\frac{\pi}{3}\}$
\item Sixth sextant: $\{\eta \in \mathbb{Z}[\omega] \mid  -\frac{\pi}{3}\leq \Arg(\eta)<0\}$
\end{enumerate}
 
The Eisenstein integers are endowed with a Euclidean function $N$ and which we will call the norm. It is defined as follows.

\begin{definition}[Norm function] $N:\mathbb{Z[\omega]} \rightarrow \mathbb{N}\cup \{0\}$ is defined by the equation $N(a+\omega b)=|a+\omega b|^2=(a+\omega b)\overline{(a+\omega b)}=a^2-ab+b^2=(a-b)^2-ab$ is the norm function in $\mathbb{Z[\omega]}$
\end{definition}

\begin{remark}
Equipped with this norm, the ring of Eisenstein integers is a Euclidean domain and thus a unique factorization domain.
\end{remark}

\begin{proposition}
$N$ is completely multiplicative.
\end{proposition}

\begin{proof}

Let $\alpha=a+\omega b$ and let $\beta=c +\omega d$. Then, $\alpha \cdot \beta=ac-bd+(ad+bc-bd)\omega$ and 

\begin{equation}
\begin{split}
    N(\alpha \cdot \beta)&=(ac-bd)^2-(ac-bd)(ad+bc-bd)+(ad+bc-bd)^2\\&=(a^2-ab+b^2)(c^2-cd+d^2)=N(\alpha)N(\beta)
\end{split}
\end{equation}
\end{proof}

\begin{proposition}
The units of $\mathbb{Z}[\omega]$ are $\pm 1, \pm \omega, \text{ and }\pm (1+\omega)$.
\end{proposition}

\begin{proof}
Suppose that $\eta$ has a multiplicative inverse. Then, $N(\eta)N(\eta^{-1})=N(\eta\eta^{-1})=N(1)=1$. Write $\eta=a+\omega b$. Then we have $N(\eta)=(a-b)^2+ab=1$. One can check that the only solutions to this equation are: $(\pm1,0)$, $(\pm1,\pm1)$, and $(0,\pm 1)$.

\end{proof}

\begin{corollary} $\epsilon$ is a unit if and only if $N(\epsilon)=1$.
\end{corollary}

For the remainder of this paper, elements of $\mathbb{Z}$ will be referred to by the name of rational integers or rational numbers and by English letters. Eisenstein integers, on the other hand, will be referred to by the name of integers or numbers and by Greek letters.

\begin{definition}[Prime]
A nonunit $\eta \in R$ is prime if, whenever $\eta | \alpha \beta$ for $\alpha,\beta \in R$, $\eta|\alpha$ or $\eta|\beta$.
\end{definition}

\begin{figure}\label{smallprimes}
\includegraphics[scale=0.3]{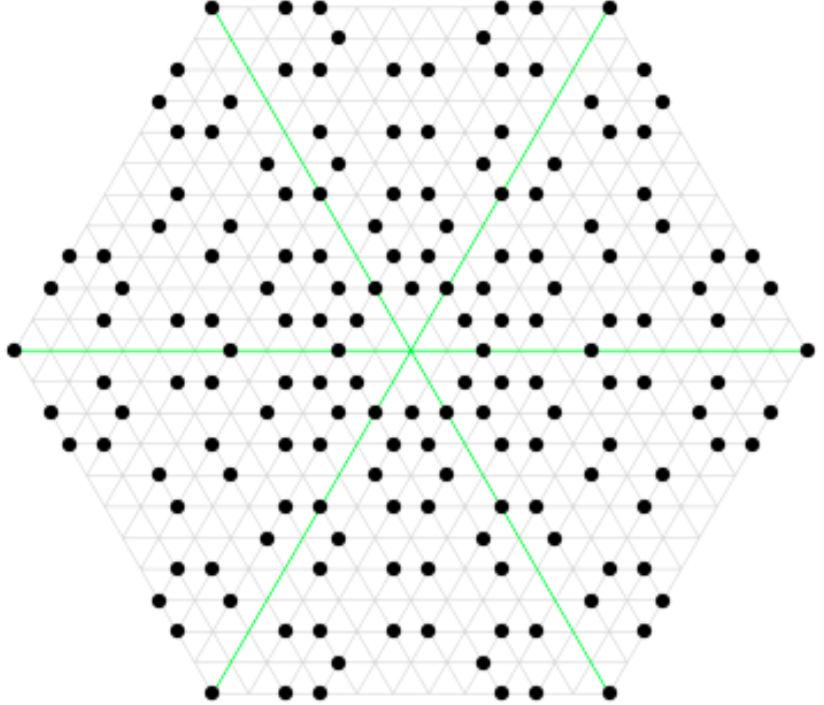}
\caption{The rays connected at the origin delimit each sextant and the black points represent Eisenstein primes.}
\label{fig:primes}
\end{figure}

For an illustration of the primes of smallest norm in the Eisenstein integers, see figure \ref{fig:primes}. The following proposition due to David Cox \cite{C} characterizes the rational primes $p$ that are also prime in the Eisenstein integers.

\begin{remark}
Remark the symmetry in figure \ref{smallprimes}. This is because if $\pi$ is prime, $\overline{\pi}$ and $\epsilon \pi$ are prime for each unit $\epsilon$.
\end{remark}

\begin{proposition} Let $p$ be a prime in $\mathbb{Z}$. Then:

\begin{enumerate}    \label{prop:primesremain}
    \item If $p=3$, then $1-\omega$ is prime in $\mathbb{Z}[\omega]$ and $3=-\omega^2(1-\omega)^2$.
    \item If $p\equiv 1 \pmod{3}$, then there is a prime $\pi \in \mathbb{Z}[\omega]$ such that $\pi \overline{\pi}$, and the primes $\pi$ and $\overline{\pi}$ are nonassociates in $\mathbb{Z}$.
    \item If $p \equiv 2 \pmod{3}$, then $p$ remains prime in $\mathbb{Z}[\omega]$.

\end{enumerate}
\end{proposition}

\begin{proposition} If $N(\eta)$ is a rational prime, then $\eta$ is prime.
\end{proposition}

\begin{proof}
Suppose that $\eta$ is not prime. Write $\eta=\alpha \beta$ for some nonunits $\alpha$ and $\beta$. $N(\eta)=N(\alpha\beta)=N(\alpha)N(\beta)$. Since $\alpha, \beta$ are nonunits, $N(\alpha), N(\beta) \geq 2$ and thus $N(\eta)$ is rational composite.
\end{proof}

\begin{definition}[Associate]
For nonzero $\eta \in R$, $\epsilon \eta$ is an associate of $\eta$ for each unit $\epsilon$. 
\end{definition}

\begin{remark}
Every nonzero $\eta$ has exactly one associate in each sextant of $\mathbb{Z}[\omega]$. 
\end{remark} 
 
 For primes $\pi \in \mathbb{Z}[\omega]$, we denote $\pi^*$ as the first-sextant associate of $\pi$. In general, for $\eta \in \mathbb{Z}[\omega]$, we define $\eta^*$ as follows. 

\begin{definition} Write $\eta=\prod_{i=1}^s \pi_i^{e_i}$ for primes $\pi$. Then, $\eta^*=\prod_{i=1}^s (\pi_i^*)^{e_i}$. 
\end{definition}

Consider, for instance, $\eta=(1-\omega)^2(\omega+3)^7$. Then, $\eta^*=(\omega+2)^2(\omega+3)^7$.

\begin{definition} [Complex sum-of-divisors function] The sum-of-divisors function \\
$\sigma: \mathbb{Z}[\omega] \rightarrow \mathbb{Z}[\omega]$ is defined by the equation

\begin{equation}
\sigma(\eta)=\sum_{\delta^*|\eta} \delta^*
\end{equation}

\end{definition}

One of the most important properties of $\sigma$ is that it is multiplicative.

\begin{remark}
If $n \in \mathbb{Z}$, then $\sigma(n)$ is the rational integers sum-of-divisors.
\end{remark}

\begin{proposition} The sum-of-divisors function is multiplicative.
\end{proposition}

\begin{proof}
Let $(\eta_1,\eta_2)=1$. We can uniquely write $\delta^*=\delta_1^*\delta_2^*$ where $\delta_1^* | \eta_1$ and $\delta_2^* | \eta_2$. Thus,

\begin{equation}
\sigma(\eta_1\eta_2)=\sum_{\delta^* | \eta_1\eta_2}\delta^*=\sum_{\delta_1^* | \eta_1,\delta_2^*|\eta_2}\delta_1^*\delta_2^*=\Bigg(\sum_{\delta_1^*|\eta_1} \delta_1^*\Bigg)\Bigg(\sum_{\delta_2^* | \eta_2} \delta_2^*\Bigg)=\sigma(\eta_1)\sigma(\eta_2)
\end{equation}

\end{proof}

\begin{definition}[$\tau$-Mersenne numbers]
For $\tau$ prime, the number

\begin{equation}
M_k=\sigma(\tau^{k-1})=\frac{\tau^k-1}{\tau-1}
\end{equation} is a $\tau$-Mersenne number. In particular, if $M_k$ is prime, it is called a $\tau$-Mersenne prime. For notational simplicity, we denote $A_k=N(M_k)$.
\end{definition}

\begin{remark} In particular, notice that if $\tau=2$, then $M_p=2^p-1$ as in the integer case.

\end{remark}

\begin{definition} Let $\eta \in \mathbb{Z}[\omega]$. $\eta$ is $\tau$-perfect if $\sigma(\eta)=\tau \eta$. $\eta$ is $\tau$-norm-perfect if $N(\sigma(\eta))=N(\tau \eta)$.
\end{definition}

\begin{remark} Every $\tau$-perfect number is norm-perfect.
\end{remark}

The following are some examples of $\tau$-norm-perfect and $\tau$-perfect Eisenstein integers for $\tau=\omega+3$. The number $\tau^{p-1}M_p$ is $\tau$-perfect for $p$ equals to 193, 709, 2029, 9049, 10453, or 255361, Clearly, for each unit $\epsilon$, we also have that $\epsilon \tau^{p-1}M_p$ is $\tau$-norm-perfect. Similarly, the number $\epsilon \tau^{p-1}\overline{M_p}$ is $\tau$-norm-perfect for $p$ equals to 11, 239, 659, 1103, and 534827.

\section{$(\omega+2)$-Perfect and $(\omega+2)$-Norm-Perfect Eisenstein Integers}\label{perfect}

In this section, we fix $\tau=\omega+2$.

Making use of the periodicity of cosine and sine, the table \ref{tab:table1} is computed.

\begin{table}[ht]
\caption{$M_k$ and $A_k$} 
\label{tab:table1}
\renewcommand\arraystretch{1.5}
\noindent\[
\begin{array}{|c|c|c|}
\hline k \pmod{12} & M_k & A_k \\ \hline
 0 & \frac{1}{2} (-1+3^\frac{k}{2})+\frac{1}{2} i (\sqrt{3}-3^{\frac{1}{2}+\frac{k}{2}}) & 1-2\cdot 3^{k/2}+3^k \\\hline
 1 & \frac{1}{2} (-1+3^{\frac{1}{2}+\frac{k}{2}})+\frac{1}{2} i (\sqrt{3}-3^\frac{k}{2}) & 1+3^k-3^{\frac{1+k}{2}}\\\hline
 2 & \frac{1}{2} (-1+2\cdot 3^\frac{k}{2})+\frac{i \sqrt{3}}{2} & 1-3^{k/2}+3^k \\\hline
 3 & \frac{1}{2} (-1+3^{\frac{1}{2}+\frac{k}{2}})+\frac{1}{2} i (\sqrt{3}+3^\frac{k}{2}) & 1+3^k \\\hline
 4 & \frac{1}{2}(-1+3^{\frac{k}{2}})+\frac{1}{2}i(\sqrt{3}+3^{\frac{1}{2}+\frac{k}{2}}) & 1+3^{k/2}+3^k \\\hline
 5 & -\frac{1}{2}+\frac{1}{2} i (\sqrt{3}+2\cdot 3^\frac{k}{2}) & 1+3^k+3^{\frac{1+k}{2}} \\\hline
 6 & \frac{1}{2} (-1-3^\frac{k}{2})+\frac{1}{2} i (\sqrt{3}+3^{\frac{1}{2}+\frac{k}{2}}) & 1+2\cdot 3^{k/2}+3^k\\\hline
 7 & \frac{1}{2} (-1-3^{\frac{1}{2}+\frac{k}{2}})+\frac{1}{2} i (\sqrt{3}+3^\frac{k}{2}) & 1+3^k+3^{\frac{1+k}{2}} \\\hline
 8 & \frac{1}{2} (-1-2\cdot 3^\frac{k}{2})+\frac{i \sqrt{3}}{2} & 1+3^{k/2}+3^k \\\hline
 9 & \frac{1}{2} (-1-3^{\frac{1}{2}+\frac{k}{2}})+\frac{1}{2} i (\sqrt{3}-3^\frac{k}{2}) & 1+3^k \\\hline
 10 & \frac{1}{2} (-1-3^\frac{k}{2})+\frac{1}{2} i (\sqrt{3}-3^{\frac{1}{2}+\frac{k}{2}}) & 1-3^{k/2}+3^k\\\hline
 11 & -\frac{1}{2}+\frac{1}{2} i (\sqrt{3}-2\cdot 3^\frac{k}{2}) & 1+3^k-3^{\frac{1+k}{2}} \\\hline
\end{array}\]

\end{table}

\begin{lemma}[Analogue of Euclid's Lemma]\label{euclid} Let $M_p$ be a Mersenne prime and $\epsilon$ a unit. If $p \equiv 1 \pmod{12}$, then $\eta=\epsilon \tau^{p-1}M_p$ is a $\tau$-norm-perfect number. If $p \equiv -1 \pmod{12}$, then $\eta=\epsilon \tau^{p-1}\overline{M_p}$ is a $\tau$-norm-perfect number.

\end{lemma}

\begin{proof} For $p \equiv 1 \pmod{12}$, $M_p$ is a sixth-sextant prime. Thus, $M_p^*=M_p(1+\omega)=\tau^p-1$. If $\eta=\epsilon \tau^{p-1}M_p$, it follows that 

\begin{equation}
N(\sigma(\eta))=N(\sigma(\epsilon)\sigma(\tau^{p-1})\sigma(M_p))=N(M_p(1+M_p^*))=N(\tau^pM_p)=N(\tau \eta)
\end{equation}

For $p \equiv -1 \pmod{12}$, $\overline{M_p}$ is a second-sextant prime. Thus, $\overline{M_p}^*=-\overline{M_p}\omega=\overline{\tau^p}-1$. If $\eta=\epsilon \tau^{p-1}\overline{M_p}$, it follows that

\begin{equation}
N(\sigma(\eta))=N(\sigma(\epsilon)\sigma(\tau^{p-1})\sigma(\overline{M_p}))=N(M_p(1+\overline{M_p}^*))=N(\overline{\tau^p}M_p)=N(\tau \eta)
\end{equation}

In both cases, $\eta$ is $\tau$-norm-perfect.

\end{proof}

Throughout the following arguments, we will make constant use of the following inequality due to McDaniel \cite{A} and improved upon by Smallbone \cite{B}.

\begin{lemma}\label{maininequality} Let $z=x+iy$ and let $k \in \mathbb{N}$.

If $x \geq \frac{5}{4}$, then

\begin{equation}
N(1+z+...+z^k)>N(z^{k-1})(N(z)+2x-1)
\end{equation}

Moreover, if $|y| \leq x-1$, then

\begin{equation}
N(1+...+z^k) \geq N(z^{k-1})(N(z)+2x+1)
\end{equation}
with equality if and only if $k=1$.

\end{lemma}

\begin{proof}

If $k=1$, 

\begin{equation}
N(1+z)=N(z)+2x+1
\end{equation}

If $k=2$,

\begin{equation}
\begin{split}
    N(1+z+z^2)&=N(z)N(z^{-1}+1+z)\\
    &=N(z)\Big(N(z)+2x+1+\frac{2x+1+2(x^2-y^2)}{N(z)}\Big)\\
    &>\begin{cases}
    N(z)(N(z)+2x-1)\text{ for all } y\\
    N(z)(N(z)+2x+1)\text{ for all } y \text{ such that } |y| \leq x
    \end{cases}
\end{split}
\end{equation}

Let $z=x+iy=re^{i \theta}$. If $k\geq3$,

\begin{equation}
\begin{split}
     N(1+z+...+z^k)&=N
     \Big(\frac{z^{k+1}-1}{z-1}\Big)=\frac{(z^{k+1}-1)(\overline{z}^{k+1}-1)}{(z-1)(\overline{z}-1)}\\
     &=\frac{N(z^{k+1})+1-(z^{k+1}+\overline{z}^{k+1})}{r^2-x+1}\\
     &=\frac{N(z^{k-1})(r^4+r^{-2(k+1)}-2r^{3-k}\cos{(k+1)\theta}}{r^2-2x+1}\\
     &>\frac{N(z^{k-1})(r^4-2)}{r^2-2x+1}
\end{split}
\end{equation}

Since $x \geq \frac{5}{4}$,

\begin{equation}
(r^2+2x-1)(r^2-2x+1)=r^4-(2x-1)^2 < r^4-2
\end{equation}

Hence, $N(1+...+z^k) > N(z^{k-1})(N(z)+2x-1)$. If also $|y| \leq x-1$, then

\begin{equation}
\begin{split}
    (r^2+2x+1)(r^2-2x+1)&=(r^2+1)^2-4x^2=r^4-2(x^2-y^2)+1\\
    & \leq r^4-(4x-3) \leq r^4-2
\end{split}
\end{equation}

Hence $N(1+z+...+z^k) > N(z^{k-1})(N(z)+2x+1)$.
\end{proof}

Many times, we will also make use of the following corollaries to lemma \ref{maininequality}.

\begin{corollary}\label{corollarytomaininequality} Let $\pi$ be prime, $k \in \mathbb{N}$, and write $\pi^*=x+iy$. Then,

\begin{equation}
\frac{N(\sigma(\pi^k))}{N(\pi^k)}>\frac{N(\pi)+2x-1}{N(\pi)}
\end{equation}

Moreover, if $y \leq x-1$, then 

\begin{equation}
\frac{N(\sigma(\pi^k))}{N(\pi^k)}\geq\frac{N(\pi)+2x+1}{N(\pi)}
\end{equation}
with equality if and only if $k=1$
\end{corollary}

\begin{corollary} For any $\eta \in \mathbb{Z[\omega]}$,

\begin{equation}\label{eq:inequalityforgeneraleta}
    \frac{N(\sigma(\eta))}{N(\eta)}\geq1
\end{equation}
with equality if and only if $\eta$ is a unit.

\end{corollary}

With these inequalities in our toolbox, we proceed lemma by lemma to prove an analogue of Euler's Lemma.

\begin{lemma}\label{firsteliminationlemma} For $k\equiv 3, 4, 5, 6, 7, 8, 9  \pmod{12}$ and $\mu$ not divisibly by $\tau$, $\eta=\tau^{k-1}\mu$ is not $\tau$-norm-perfect.
\end{lemma}

\begin{proof} Consulting table \ref{tab:table1}, it follows that, for $k \equiv 3, 4, 5, 6, 7, 8, 9 \pmod{12}$, $A_k>3^k=N(\tau^k)$. Thus, by inequality (\ref{eq:inequalityforgeneraleta}), it follows that

\begin{equation}
N(\sigma(\eta))=A_kN(\sigma(\mu))>N(\tau^k)N(\mu)=N(\tau \eta)
\end{equation}

Hence, $\eta$ is not $\tau$-norm-perfect.

\end{proof}

We summarize the results of lemma \ref{firsteliminationlemma} in the following corollary.

\begin{corollary} \label{corollaryelimination}If $\eta=\tau^{k-1} \mu$ is $\tau$-norm-perfect, then $k \equiv 0, \pm 1, \pm 2 \pmod{12}$.
\end{corollary}

\begin{lemma}\label{Mkisprimelemma} Let $k\geq 2$ and $\mu$ not divisible by $\tau$. If $\eta=\tau^{k-1}\mu$ is $\tau$-norm-perfect, then $M_k$ or $\overline{M_k}$ divide $\eta$ and are both prime.
\end{lemma}

\begin{proof}

Let $\pi$ be a first-sextant prime divisor of $M_k$. Suppose that $\eta$ is $\tau$-norm-perfect. Then, it follows that

\begin{equation}
    \begin{split}
        3 n \overline{\eta}=N(\tau \eta)=N(\sigma(\eta))=N( M_k \sigma( \mu ))=\pi \overline{\pi} N\Big( \frac{M_k}{\pi} \sigma( \mu )\Big )
    \end{split}
\end{equation}

Thus, it follows that $\pi| 3 n \overline{\eta}$. Since $3=(1+w)(1-w)^2$, since $1+w$ is a unit, since $1-\omega$ is an associate of $\tau$ and since $(M_k,\tau)=1$, it follows that $\pi \nmid 3$. Thus,  $\pi | n \overline{\eta}$. Since $\pi$ is prime, then it follows that $\pi | \eta$ or $\pi | \overline{\eta}$. Equivalently, $\pi | \eta$ or $\overline{\pi} | \eta$. In particular, since $(M_k,\tau)=1$, it follows that $\pi | \mu$ or $\overline{ \pi } | \mu$.

For any prime $\pi$ such that $\pi | \mu$, let $a$ be the largest rational integer such that $\pi^a | \mu$. Using corollaries to lemma \ref{maininequality}, it follows that

\begin{equation}
    \begin{split}
        1=\frac{N(\sigma(\eta))}{N(\tau \eta)}&=\frac{N(\sigma(\tau^{k-1})\sigma(\pi^a))}{N(\tau^k\pi^a)}\frac{N(\sigma(\mu/\pi^a))}{N(\mu/\pi^a)}\geq \frac{N(\sigma(\tau^{k-1})\sigma(\pi^a))}{N(\tau^k\pi^a)}\\
        &>\frac{A_k(N(\pi)+2x-1)}{N(\tau^k)N(\pi)}
    \end{split}
\end{equation}

Rearranging gives us

\begin{equation}
N(\pi)>\frac{A_k(2x-1)}{N(\tau^k)-A_k}
\end{equation}

Since $\pi^*$ and $ \overline{\pi}^*$ are first-sextant primes different from $\tau$, $\R \pi^*, \R \overline{\pi}^* \geq 2$. So

\begin{equation}
N(\pi)>\frac{3A_k}{N(\tau^k)-A_k}
\end{equation}

By corollary \ref{corollaryelimination}, it follows that

\begin{equation}
\frac{3A_k^{\frac{1}{2}}}{N(\tau^k)-A_k} \geq \frac{3 \left(3^{\frac{k}{2}}-1\right)}{2\cdot 3^{\frac{k}{2}}-1}>1
\end{equation}

Thus,

\begin{equation}
N(\pi) > \frac{3A_k}{N(\tau^k)-A_k} > A_k^{\frac{1}{2}}=N(M_k)^{\frac{1}{2}}
\end{equation}

That is, 

\begin{equation}
N(\pi)^2>N(M_k)
\end{equation}

Assume that $M_k$ is not prime. Write $M_k=\epsilon \pi_0 \pi_1...\pi_r$ for $r \in \mathbb{N}$, where $\pi_i$ is a prime, and $\epsilon$ is a unit. Let $\pi_0$ be a prime with the least norm among the norm of the primes $\pi_i$. Then, it follows that

\begin{equation}
N(\pi_0)>N(\pi_1)...N(\pi_r)
\end{equation}which is a contradiction. Thus, it follows that $M_k=\epsilon \pi$ for some prime $\pi$ and unit $\epsilon$.

Suppose that $\overline{M_k}$ is not prime. Write $\overline{M_k}=\alpha \beta$. Then, $M_k=\overline{\alpha} \overline{\beta}$, making $M_k$ not prime. Therefore, by the above argument, $M_k$ and $\overline{M_k}$ are both prime.

\end{proof}

\begin{lemma}\label{primepower}
 If $M_k$ is prime, then $k$ is rational prime.
\end{lemma}

\begin{proof}Suppose that $k$ is composite. Write $k=nm$ for $n,m \geq 2$. Then,

\begin{equation}
M_k=\frac{\tau^k-1}{\tau-1}=\frac{\tau^{nm}-1}{\tau-1}=\Bigg(\frac{\tau^{n}-1}{\tau-1}\Bigg)\Bigg(\frac{\tau^{nm}-1}{\tau^n-1}\Bigg)
\end{equation}

If $\frac{\tau^{nm}-1}{\tau^n-1}=\epsilon$ for some unit $\epsilon$, then by rearranging and taking norms, it follows that

\begin{equation}
3^nN(1-\epsilon \tau^{mn-n})=N(1-\epsilon)
\end{equation}but $3^nN(1-\epsilon \tau^{mn-n})\geq9$ and $N(1-\epsilon) \leq 4$. By the same argument, $\frac{\tau^{n}-1}{\tau-1}$ is not a unit.

\end{proof}

\begin{lemma}\label{lemmaonformofnormperfect} Let $t \in \mathbb{N}$, $\delta$ not divisibly $\tau$, and $k \geq 2$. If $\eta=\tau^{k-1}\mu$ is a $\tau$-norm-perfect number, then, for some unit $\epsilon$, either $\eta=\epsilon \tau^{p-1} M_p^t\delta$ where $M_p$ is a Mersenne prime with $p \equiv 1 \pmod{12}$, or $\eta=\epsilon \tau^{p-1} \overline{M_p}^t\delta$ where $M_p$ is a Mersenne prime with $p \equiv -1 \pmod{12}$
\end{lemma}

\begin{proof} By lemma \ref{Mkisprimelemma}, $\eta=\tau^{k-1}M_k^t \delta$ or $\eta=\tau^{k-1}\overline{M_k}^t \delta$ for some $\delta$ not divisible by $\tau$. By choosing $t$ sufficiently large, we get that $( M_k , \delta)=1$ or $( \overline{M_k} , \delta)=1$, respectively. By proposition \ref{primepower}, $k$ must be a rational prime. Hence, we write $p$. By corollary \ref{corollaryelimination}, $p=2$ or $p \equiv \pm 1 \pmod{12}$.

We are left to show that for $p=2$ and $M_p$ prime, $\eta=\epsilon \tau^{p-1} M_p^t\delta$ and $\eta=\epsilon \tau^{p-1} \overline{M_p}^t\delta$ are not $\tau$-norm-perfect; that, for $p \equiv -1 \pmod{12}$ and $M_p$ prime, $\eta=\epsilon \tau^{p-1} M_p^t\delta$ is not $\tau$-norm-perfect; and that, for $p \equiv 1 \pmod{12}$ and $M_p$ prime, $\eta=\epsilon \tau^{p-1} \overline{M_p}^t\delta$ is not $\tau$-norm-perfect.

Consider $\eta=\tau M_2^t \delta$. $M_2=\sigma(\tau)=1+\tau=3+\omega$. $M_2^*=M_2$. So, by lemma \ref{maininequality} and its corollary, it follows that

\begin{equation}
    \begin{split}
        \frac{N(\sigma(\eta))}{N(\tau \eta)}&= \frac{N(\sigma(\tau))}{N(\tau)^2}\frac{N(\sigma(M_2^t))}{N(M_2^t)}\frac{N(\sigma(\delta))}{N(\delta)}\geq\frac{N(\sigma(\tau))}{N(\tau)^2}\frac{N(\sigma(M_2^t))}{N(M_2^t)}\\
        &>\frac{N(1+\tau)}{N(\tau)^2}\frac{A_2+2\R M_2^*-1}{A_2}=\frac{11}{9}>1
    \end{split}
\end{equation}

Consider $\eta=\tau \overline{M_2}^t \delta$. $\overline{M_2}=\overline{\sigma(\tau)}=\overline{1+\tau}=\overline{3+\omega}=2-\omega$. $\overline{M_2}^*=\overline{M_2}(\omega +1)=3+2\omega$. As before, it follows that

\begin{equation}
    \begin{split}
        \frac{N(\sigma(\eta))}{N(\tau \eta)}>\frac{N(1+\tau)}{N(\tau)^2}\frac{A_2+2\R \overline{M_2}^*-1}{A_2}=\frac{10}{9}>1
    \end{split}
\end{equation}

Consider $\eta=\epsilon \tau^{p-1} M_p^t\delta$ for $p\equiv -1 \pmod{12}$. Since $M_p$ is a fifth-sextant prime, $M_p^*=\omega M_p$. Since $\I M_p^* \leq \R M_p^* -1$, it follows that

\begin{equation}
    \begin{split}
        \frac{N(\sigma(\eta))}{N(\tau \eta)}\geq\frac{A_p+2 \R M_p^*+1}{N(\tau^p)}=3^p-1>1
    \end{split}
\end{equation}

Consider $\eta=\epsilon \tau^{p-1} \overline{M_p}^t\delta$ for $p\equiv 1 \pmod{12}$. Since $M_p$ is a sixth-sextant prime, $\overline{M_p}^*=\overline{M_p}$. Since $\I \overline{M_p}^* \leq \R \overline{M_p}^* -1$, it follows that

\begin{equation}
    \begin{split}
        \frac{N(\sigma(\eta))}{N(\tau \eta)}\geq\frac{A_p+2 \R \overline{M_p}^*+1}{N(\tau^p)}=3^p-1>1
    \end{split}
\end{equation}

\end{proof}

We now present the analogue of Euler's lemma.

\begin{lemma}[Analogue of Euler's Lemma] Let $k \geq 2$. If $\eta=\tau^{k-1}\eta$ is a $\tau$-norm-perfect number, then, for some unit $\epsilon$, either $\eta=\epsilon \tau^{p-1} M_p$ where $M_p$ is a Mersenne prime with $p \equiv 1 \pmod{12}$, or $\eta=\epsilon \tau^{p-1} \overline{M_p}$ where $M_p$ is a Mersenne prime with $p \equiv -1 \pmod{12}$.
\end{lemma}

\begin{proof}
Let $M_p$ prime and $p \equiv 1 \pmod{12}$. Since $| \I M_p^*| \leq \R M_p^* -1$, by corollary \ref{corollarytomaininequality},

\begin{equation}\label{Mpresult}
\frac{N(\sigma(M_p^t))}{N(M_p^t)} \geq\frac{N(\sigma(M_p))}{N(M_p)}
\end{equation}
with equality if and only if $t=1$.

By the same argument we also have that 

\begin{equation}
\frac{N(\sigma(\overline{M_p^t}))}{N(\overline{M_p^t})} \geq \frac{N(\sigma(\overline{M_p}))}{N(\overline{M_p})}
\end{equation}
with equality if and only if $t=1$.

Suppose that $\eta$ is $\tau$-norm-perfect number, then, by lemma \ref{lemmaonformofnormperfect}, $\eta=\epsilon \tau^{p-1} M_p^t \delta$ or $\eta=\epsilon \tau^{p-1} \overline{M_p}^t \delta$. 

Assume that $\eta$ is of the former form. Then, by the Analogue of Euclid's Lemma, by corollary \ref{corollarytomaininequality}, by inequality \ref{Mpresult}, and since $\eta$ is $\tau$-norm-perfect, it follows that

\begin{equation}
\begin{split}
1&=\frac{N(\sigma(\eta))}{N(\tau \eta)}=\frac{N(\sigma(\tau^{p-1}))}{N(\tau^{p})}\frac{N(\sigma(M_p^t))}{N(M_p^t)}\frac{N(\sigma(\delta))}{N(\delta)}\\ &\geq\frac{N(\sigma(\tau^{p-1}))}{N(\tau^{p})}\frac{N(\sigma(M_p))}{N(M_p)}\frac{N(\sigma(\delta))}{N(\delta)}\\
&=\frac{N(\sigma(\tau^{p-1}M_p))}{N(\tau^{p}M_p)}\frac{N(\sigma(\delta))}{N(\delta)}=\frac{N(\sigma(\delta))}{N(\delta)}
\end{split}
\end{equation}

Thus, $\frac{N(\sigma(\delta))}{N(\delta)}=1$; that is, $\delta$ is a unit. Further, if $\delta$ is a unit, it follows that $\frac{N(\sigma(M_p^t))}{N(M_p^t)}= \frac{N(\sigma(M_p))}{N(M_p)}$; that is, that $t=1$. By the same argument, it follows that $\delta$ is a unit and $t=1$ in the latter form of $\eta$.

\end{proof}

We consolidate the analogues of Euclid's and Euler's lemmas into what we have called the Euclid-Euler Theorem for $\tau$-norm-perfect Eisenstein Integers.

\begin{theorem}[Euclid-Euler Theorem for $\tau$-Norm-Perfect Eisenstein Integers]\label{normperfect} Let $M_p$ be a Mersenne prime and $\epsilon$ a unit. If $ p\equiv 1 \pmod{12}$, $\eta=\epsilon \tau^{p-1} M_p$ is a $\tau$-norm-perfect number; if $p \equiv -1 \pmod{12}$, $\eta=\epsilon \tau^{p-1} \overline{M_p}$ is a $\tau$-norm-perfect number. Conversely, if $\eta$ is a $\tau$-norm-perfect number divisible by $\tau$, then, for some unit $\epsilon$, either $\eta=\epsilon \tau^{p-1} M_p$, where $M_p$ is a Mersenne prime with $p \equiv 1 \pmod{12}$, or $\eta=\epsilon \tau^{p-1} \overline{M_p}$, where $M_p$ is a Mersenne prime with $p \equiv -1 \pmod{12}$.
\end{theorem}

\begin{corollary} There are no imprimitive $\tau$-Norm-Perfect numbers divisible by $\tau$ in the Eisenstein integers.

\end{corollary}

We now derive what we have called the Euclid-Euler Theorem for $\tau$-Perfect Eisenstein Integers.

\begin{corollary}[Euclid-Euler Theorem for $\tau$-Perfect Eisenstein Integers] Let $M_p$ be a Mersenne prime. Then, $\eta$ is an $\tau$-perfect number divisible by $\tau$ if and only if $\eta=\tau^{p-1} M_p$ for $p \equiv 1 \pmod{12}$.
\end{corollary}

\begin{proof} Consider $\eta=\epsilon \tau^{p-1}M_p$ for $p \equiv 1 \pmod{12}$. Since $M_p$ is a sixth-sextant prime, $M_p^*=M_p(1+\omega)=\tau^p-1$. Thus, it follows that 

\begin{equation}
\sigma(\eta)=\sigma(\tau^{p-1}M_p)=\sigma(\tau^{p-1})\sigma(M_p)=M_p(1+M_p^*)=\tau^pM_p=\tau \eta
\end{equation}

By theorem \ref{normperfect}, if $\eta$ is an $\tau$-perfect number divisible by $\tau$, then $\eta=\epsilon \tau^{p-1} M_p$ for $p \equiv 1 \pmod{12}$, $M_p$ prime, and some unit $\epsilon$; or  $\eta=\epsilon \tau^{p-1} \overline{M_p}$ for $p \equiv -1 \pmod{12}$, $M_p$ prime, and some unit $\epsilon$. 

Consider the latter. Since, for $p \equiv -1 \pmod{12}$, $M_p$ is a fifth-sextant prime, $\overline{M_p}^*=-\omega \overline{M_p}=\overline{\tau^p}-1$. Thus, it follows that

\begin{equation}
\sigma(\eta)=M_p(1+M_p^*)=\overline{\tau^p}M_p
\end{equation}
         
Since $\tau^p \neq \overline{\tau^p}$, $\eta$ is not $\tau$-perfect. Therefore, if $\eta$ is $\tau$-perfect, then $\eta=\epsilon \tau^{p-1} M_p$ for $p \equiv 1 \pmod{12}$, $M_p$ prime, and some unit $\epsilon$. It is easy to check that $\eta$ is only $\tau$-perfect for $\epsilon=1$.

\end{proof}

\section{Nonexistence of $2$-norm-perfect Eisenstein Integers}\label{smallbonesresults}

\begin{lemma} Let $k\geq 2$. If $\eta=2^{k-1}\mu$ is $2$-norm-perfect, then $\sigma(2^{k-1})$ is prime.
\end{lemma}

\begin{proof}
Let $\pi$ be a first-sextant prime factor of $\sigma(2^{k-1})$ with $\R \pi \geq 2$. Let $a$ be the largest rational integer such that $\pi^a | \sigma(2^{k-1})$. 

\begin{equation}
2^2 \eta \overline{\eta}=N(2\eta)=N(\sigma(2^{k-1})\sigma(\mu))=\pi \overline{\pi} N\Big(\frac{\sigma(2^{k-1})}{\pi}\sigma(\mu)\Big)
\end{equation}

Thus it follows that $\pi | 2^2 \eta \overline{\eta}$. Since $(\pi,2)=1$, it follows that $\pi | \eta \overline{\eta}$. Since $\pi$ is prime. $\pi | \eta$ or $\overline{\pi} | \eta$. As shown in the proof of lemma \ref{Mkisprimelemma},

\begin{equation}
N(\pi)>\frac{3(2^k-1)^2}{2^{2k}-(2^k-1)^2}=\frac{3(2^k-1)^2}{2^{k+1}-1}
\end{equation}

Since

\begin{equation}
\frac{3(2^k-1)}{2^{k+1}-1}>1
\end{equation}

it follows that 

\begin{equation}
N(\pi)>\sigma(p^{k-1})
\end{equation}

If $\sigma(p^{k-1})$ is composite and write $\sigma(p^{k-1})=\pi_0\pi_1...\pi_r$ where $\pi_0$ is a prime with the least norm among the norm of all prime factors of $\sigma(p^{k-1})$, then

\begin{equation}
N(\pi_0)>N(\pi_1)...N(\pi_r)
\end{equation}

This is a contradiction. Thus, $\sigma(p^{k-1})$ is prime.
\end{proof}

\begin{theorem}
There are no $2$-norm-perfect Eisenstein integers divisible by $2$.
\end{theorem}

\begin{proof}
By proposition \ref{prop:primesremain}, $\sigma(2^{k-1})$ is not prime since $\sigma(2^{k-1})=1+2+...+2^{k-1}=2^k-1 \not \equiv 2 \pmod{3}$.
\end{proof}

\section{Discussion}\label{discussion}

In regard to future work, we are interested in studying $\tau$-norm-perfect and $\tau$-perfect numbers for other values of $\tau$, and in studying $\tau$-norm-perfect and $\tau$-perfect numbers that are not divisible by $\tau$ for $\tau=2$ and $\tau=\omega+2$. Thus far, in the Gaussian and in the Eisenstein integers, there are only characterizations for $\tau$-norm-perfect and $\tau$-perfect numbers that are also divisible by $\tau$. 

In the Eisenstein integers, two potentially promising start points are to attempt to characterize or prove nonexistence of $\tau$-norm-perfect integers for $\tau=\omega+3$ or $2\omega+3$ as these are the first-sextant primes that follow $\omega+2$ and $2$ in norm. Alternaively, rational primes $p \equiv 2 \pmod{3}$ may be studied. We also would like to remark that many of our computations become infeasible after certain modifitication. In particular, computing $\sigma(\eta^k)$ for nonprime $\eta$ and general $k \in \mathbb{N}$ usually turns to be a cumbersome task.

\bibliographystyle{amsplain}

\end{document}